\documentclass[a4paper,11pt]{article}
\usepackage{latexsym,amsmath}

\author{\emph{\textbf{Shaohua Zhang}}$^{1, 2}$\\\\
1 School of Mathematics, Shandong University,\\
Jinan,  Shandong, 250100, China \\
2 The key lab of cryptography technology and information
security, \\
Ministry of Education, Shandong University, \\
Jinan, Shandong, 250100, China\\
E-mail: shaohuazhang@mail.sdu.edu.cn}
\title{\textbf{On Fixed Points of Order K of RSA}}
\date{}

\begin{document}
\maketitle

\begin{abstract}
In this paper, we gave a preliminary dynamical analysis on the RSA
cryptosystem and obtained a computational formulae of the number of
the fixed points of $k$ order of the RSA. Thus, the problem in [8,
9] has been solved.

\vspace{3mm}\textbf{Keywords:}  RSA; fixed point; fixed points of
order $k$;  fixed points attack; dynamical analysis

\vspace{3mm}\textbf{2000 MR  Subject Classification:}\quad 11A25;
11T71; 37A45
\end{abstract}

\section{Introduction}
\setcounter{section}{1}\setcounter{equation}{0}

Shortly after Diffie and Hellman [1] introduced the idea of public
key cryptography, Rivest, Shamir and Adleman (RSA) [2] proposed such
a cryptosystem. A simplified version of RSA is the following:

\vspace{3mm} Let $n=pq$ be the product of two large primes of the
same size. Let $e,d$ be two integers satisfying $ed\equiv 1(\mod
\varphi (n))$. Call $n$ the RSA modulus, $e$ the encryption
exponent, and $d$ the decryption exponent. Let $e$ and $n$ be public
keys, and let $d$ be the corresponding secret key. A message is an
integer $m\in Z_{n}$. To encrypt $m$, one computes $m^{e}\equiv c(
\mod n)$. To decrypt the ciphertext $c$, the receiver computes
$m\equiv m^{ed}\equiv c^{d}(\mod n)$. Denote such a cryptosystem by
$RSA(n,e)$. We call $m$ a fixed point of $RSA(n,e)$ if $m^{e}\equiv
m(\mod n)$. And call $m$ a fixed point of order $k$ if $k$ is the
smallest positive integer such that $m^{e^{k}}\equiv m(\mod n)$.
Clearly, $f:x\longrightarrow x^{e}(\mod n)$ is a dynamical system.
Thus, $k$ is exactly the period of $m$. For more details on the
arithmetic of dynamical systems, see [7].

\vspace{3mm} In 1979, Blakley and Borosh [3] first pointed out that
there were at least 9 fixed points in $RSA(n,e)$. For more
references on fixed points, also see [4]-[6]. Denoted the set of all
fixed points of order $k$ of $RSA(n,e)$ by $E_{n,e,k}$ and the
cardinality of the set $S$ by $|S|$. In [8, 9], Yu considered the
general case of fixed points of order $k$ and gave geometric mean
value of $|T_{n,e,k}|$ and pointed out that it was difficult to give
a quantitative description of $|T_{n,e,k}|$, where $k$ is a given
positive integer and
$$T_{n,e,k}=\{x|\forall m<k,m\in N,x\in Z_{n}^{\ast },x^{e^{k}}\equiv
x(\mod n),x^{e^{m}}\neq x(\mod n)\}.$$

\vspace{3mm}In this essay, we preliminarily consider this question
and obtain the following results:

\vspace{3mm}\noindent{\bf Theorem ~1~~}%
$|T_{n,e,k}|=\sum_{d|k}\mu (k/d)(e^{d}-1,p-1)(e^{d}-1,q-1)$ , where
$\mu (\cdot )$ is the M\"{o}bius function.

\vspace{3mm}Based on this result, we get Theorem 2.

\vspace{3mm}\noindent{\bf Theorem ~2~~}%
$|E_{n,e,k}|=\sum_{d|k}\mu (k/d)((e^{d}-1,p-1)+1)((e^{d}-1,q-1)+1)$.

\section{Proof of Main Theorems}
We denote the set of positive integers by $N$ . For given positive
integers $a$ and $b$, we write $a|b$ if $a$ divides $b$ . And denote
the greatest common divisor of $a$ and $b$ by $(a,b)$. Denote a
complete set of residues modulo $n$ by $Z_{n}$, where $1<n\in N$ ,
and a reduced set of residues modulo $n$ is denoted by $Z_{n}^{\ast
}$. Let $a$ be an integer relatively prime to $n$. The order of $a$
modulo $n$, denoted by $ord_{n}(a)$, which is the smallest positive
integer $d$ such that $a^{d}\equiv 1(\mod n)$.

\vspace{3mm}\noindent{\bf Lemma 1[8]~~} %
For $1<n\in N$, $r\in N$, let the canonical factorization of $n$ be
$\prod\limits_{i=1}^{m}p_{i}^{a_{i}}$ and $T_{n,r}=\{x|x^{r}\equiv
1( \mod n),1\leq x<n\}$, then
$|T_{n,r}|=\prod\limits_{i=1}^{m}(r,\varphi (p_{i}^{a_{i}}))$.

\vspace{3mm}\noindent{\bf Lemma~2~~}%
For $1<n\in N$, $a,m,k\in N$, $e\in Z_{\varphi (n)}^{\ast }$, if
$a^{e^{k}}\equiv a(\mod n)$ and $k|m$, then $a^{e^{m}}\equiv a(\mod
n)$.

\vspace{3mm}\noindent{\bf Proof~~}%
Let $m=tk$. When $t=1$, clearly $a^{e^{k}}\equiv a^{e^{m}}\equiv
a(\mod n)$. Suppose that $a^{e^{m}}\equiv a(\mod n)$ when $t=l$. And
when $t=l+1$, we have $a^{e^{m}}\equiv a^{e^{lk}e^{k}}\equiv
a^{e^{k}}\equiv a(\mod n)$. It immediately shows that Lemma 2 is
true by induction.

\vspace{3mm}\noindent{\bf Lemma~3~~}%
For $1<n\in N$, $a,m,k\in N$, $e\in Z_{\varphi (n)}^{\ast }$, if
$a^{e^{m}}\equiv a(\mod n)$ and $a\in E_{n,e,k}$, then $k|m$.

\vspace{3mm}\noindent{\bf Proof ~~}%
Let $m=kt+r$, $t\in N$, $0\leq r<k$. We have $a^{e^{m}}\equiv
a^{e^{kt}e^{r}}\equiv a^{e^{r}}\equiv a(\mod n)$ by Lemma 2. Since
$a\in E_{n,e,k}$, hence $r=0$, and Lemma 3 is true.

\vspace{3mm}\noindent{\bf Proof of Theorem 1~~}%
By Lemma 1 and Lemma 3, it is easy to deduce
$\sum_{d|k}|T_{n,e,d}|=\prod\limits_{i=1}^{m}(e^{k}-1,\varphi
(p_{i}^{a_{i}}))$. By M\"{o}bius inversion, it immediately shows
that Theorem 1 is true.

\vspace{3mm}\noindent{\bf Proof of Theorem 2~~}%
By Lemma 2 and Lemma 3, analogously, using Chinese Remainder Theorem
and the method of proof of Theorem 1, it is easy to deduce that
Theorem 2 is true.

\vspace{3mm}\noindent{\bf Corollary~1~~}%
Let $1<n\in N$, $r\in N$, and let the canonical factorization of $n$
be $\prod\limits_{i=1}^{m}p_{i}^{a_{i}}$. Then
$|\{x|ord_{n}(x)=r,1\leq x\in Z_{n}^{\ast }\}|=\sum_{d|r}(\mu
(r/d)\prod\limits_{i=1}^{m}(d,\varphi (p_{i}^{a_{i}})))$.

\vspace{3mm}\noindent{\bf Corollary~2~~}%
Let $1<n\in N$, $r\in N$, let the canonical factorization of $n $ be
$\prod\limits_{i=1}^{m}p_{i}^{a_{i}}$, and let $F_{n,r}=\{x|\forall
k<r,k\in N,1\leq x\leq n,x^{r}\equiv x(\mod n),x^{k}\neq x(\mod
n)\}$. Then $|F_{n,r}|=\sum_{d|r}(\mu
(r/d)\prod\limits_{i=1}^{m}(1+(d-1,\varphi (p_{i}^{a_{i}}))))$.

\section{Conclusion}
Clearly, if the factorization of $n$ is known, then computing the
number of the fixed points of order $k$ of the RSA cryptosystem is
simple and convenient by the presented formulae. This is useful to
pick the encryption exponent, which is necessary to ensure the
resulting RSA safe from fixed points attack. Maybe we are not afraid
of a fixed point. However, the following problem should be further
considered: Is there a polynomial-time algorithm for finding a fixed
point $m$, where $m\neq 0,\pm1$? This problem and Factoring the RSA
modulus perhaps are equivalent.

\vspace{3mm}\noindent\textbf{Remark:} This paper is the revision of
paper [10] in the proceedings of China Crypt'2006, whose Chinese
version has been accepted by \emph{Journal of Mathematics }(Wuhan,
China).

\section{Acknowledgements}
I am thankful to the referees for their suggestions improving the
presentation of the paper and also to my supervisor Professor Wang
Xiaoyun for her valuable help and encouragement. Thank Institute for
Advanced Study in Tsinghua University for providing us with
excellent conditions. This work was partially supported by the
National Basic Research Program (973) of China (No. 2007CB807902)
and the Natural Science Foundation of Shandong Province (No.
Y2008G23).

\clearpage

\end{document}